\documentclass[a4paper, conference]{IEEEtran}
\IEEEoverridecommandlockouts
\usepackage{amsfonts}
\usepackage[cmex10]{amsmath}
\usepackage{graphicx}
\usepackage{setspace} 
\usepackage{times} 
\usepackage{epsfig}
\usepackage{subfigure} 
\usepackage{latexsym} 
\usepackage{cite}
\usepackage{epsfig}
\usepackage{fixltx2e}
\usepackage{verbatim}

\ifCLASSINFOpdf
\else
\fi

\begin{document}
\title{\huge Adaptive Sensing and Transmission Durations for\\ Cognitive Radios}
\author{\IEEEauthorblockN{Wessam Afifi, Ahmed Sultan and Mohammed Nafie}
\thanks{$^{1}$This work was supported in part by a grant from the Egyptian NTRA (National Telecommunications Regulatory Authority).}
\IEEEauthorblockA{Wireless Intelligent Networks Center (WINC)\\
Nile University, Cairo, Egypt.\\
E-mail: wessam.afifi@nileu.edu.eg, \{asultan, mnafie\}@nileuniversity.edu.eg }}
\maketitle
%----------------------------------------------------Abstract-------------------------------------------------------------------
\begin{abstract}
In a cognitive radio setting, secondary users opportunistically access the spectrum allocated to primary users. Finding the optimal sensing and transmission durations for the secondary users becomes crucial in order to maximize the secondary throughput while protecting the primary users from interference and service disruption. In this paper an adaptive sensing and transmission scheme for cognitive radios is proposed. We consider a channel allocated to a primary user which operates in an unslotted manner switching activity at random times. A secondary transmitter adapts its sensing and transmission durations according to its belief regarding the primary user state of activity. The objective is to maximize a secondary utility function. This function has a penalty term for collisions with primary transmission. It accounts for the reliability-throughput tradeoff by explicitly incorporating the impact of sensing duration on secondary throughput and primary activity detection reliability. It also accounts for throughput reduction that results from data overhead. Numerical simulations of the system performance demonstrate the effectiveness of adaptive sensing and transmission scheme over non-adaptive approach in increasing the secondary user utility.$^{1}$
\end{abstract}
\IEEEpeerreviewmaketitle
%----------------------------------------------------Introduction---------------------------------------------------------------
\section{Introduction}
In cognitive radio networks, secondary or unlicensed users are allowed to share the spectrum with primary, licensed users. Unlike primary users who can access the spectrum at will at any time, secondary users have to search for the vacant slots in the spectrum and opportunistically access the spectrum without causing interference to the primary users. Secondary users carry out spectrum sensing to detect the state of the primary users. When there is no primary activity over a certain band, the secondary users can utilize the band to transmit their own data. 

One important aspect in cognitive radio networks that has received wide attention in research is to find the optimal sensing and transmission strategies for the secondary user. This includes, inter alia, the determination of the optimal inter-sensing time for unslotted primary networks \cite{kim}, \cite{intersensing_duration}, specifying the optimal channel order for sensing and access \cite{sensing_order_michigan}, \cite{optimal_channel_order_Poor}, and finding the optimal sensing duration based on secondary observations \cite{Opportunistic_spectrum_access}. Note that the detection of primary users becomes more reliable as the sensing duration increases. On the other hand, and under the assumption that a secondary user either senses or transmits over a channel, a long sensing duration means decreasing the time available for transmission. Sensing constitutes an overhead that comes at the expense of transmission \cite{sensing_overhead}. This is true in both slotted and unslotted primary systems \cite{kim}, \cite{Opportunistic_spectrum_access} indicating the existence of a tradeoff between sensing reliability and secondary throughput \cite{sensing_throughput_tradeoff}. 

The authors of \cite{senhua_1} develop their preliminary work in \cite{senhua_2} and design a secondary access scheme that optimizes the secondary access efficiency while protecting the primary transmission from interference. The primary mode of operation is un-slotted, which means that its active and idle times are random variables. A utility function is developed to account for the secondary throughput and to penalize it for colliding with primary transmission. During the idle primary period, the secondary can either sense or transmit. Both the sensing and transmission durations are assumed to be fixed. The optimal solution is threshold-based such that the secondary transmits when its belief about the primary being idle exceeds a certain threshold. The belief is updated based on the secondary sensing observations and the feedback it receives from its respective receiver. The authors assume that the secondary user can perfectly detect the start point of the primary off duration. The problem of the quickest detection of transmission opportunity is addressed in \cite{Zhao_Quickest_Change_Detection} and \cite{zhao_3}, and is beyond the scope of this work.

We build on the work in \cite{senhua_1} and make the following contributions. Instead of using fixed sensing and transmission durations for the secondary users, we consider varying the durations according to the belief of the secondary user concerning the primary activity. The durations become optimization variables that parameterize secondary utility function. The motivation for this is that the secondary transmitter may waste time and energy in long sensing periods although it has a high belief that the primary user is idle. On the other hand, the secondary transmitter may relatively increase its sensing period to detect the actual state of the primary user if it has a considerable belief that the primary user is busy. This long sensing duration makes the sensing outcome more reliable and reduces the probability of collision with the primary. Therefore, adaptive sensing and transmission durations can enhance the secondary throughput and afford more protection to the primary compared to the case of fixed durations. In addition, previous works also consider the secondary throughout on the basis of the whole transmission duration. This ignores the data overhead which makes it better for the secondary user to make one relatively long transmission instead of making multiple small transmissions. In a cognitive setting, however, a long transmission duration increases the probability of colliding with the primary user. Thus there is another tradeoff here if the overhead is accounted for as we do in this paper. We also allow the secondary to remain idle because in practice both sensing and transmission have a cost in terms of power consumption. 

The rest of the paper is organized as follows. The system model and problem formulation are described in Section \ref{section:system_model}. We present the adaptive sensing and transmission scheme in Section \ref{section:adaptive_structure}. In Section \ref{section:simulation_results} we provide simulation results and compare between the non-adaptive and adaptive schemes. We conclude the paper in Section \ref{section:conclusion}. 
%---------------------------------------------------System Model----------------------------------------------------------------
\section{System Model and Problem Formulation}{\label{section:system_model}}
We consider a channel allocated to a primary Tx-Rx pair which operates in an unslotted manner, switching activity at random times. A secondary terminal attempts to opportunistically access this channel maximizing its throughput while simultaneously minimizing the probability of colliding with primary transmission.

The primary user's activity follows an alternating on/off renewal process with certain probability distributions for the idle and busy periods: $f_{\rm off}\left(t\right)$ and $f_{\rm on}\left(t\right)$ with means $T_{\rm off}$ and $T_{\rm on}$, respectively. The idle and busy periods are independent of each other. We assume that there is no cooperation between the primary and secondary users. The secondary user can quickly and reliably detect the transition of the primary user from busy to idle. This transition represents $t=0$. 

\subsection{Secondary Actions}
We assume that the secondary user always has data to transmit. Its objective is to enhance its transmission throughput while protecting the primary user from interference. The secondary transmitter can perform one of three actions: stay idle $\left\{I\right\}$, carry out spectrum sensing $\left\{S\right\}$, or transmit its data $\left\{T\right\}$ with action space $A=\left\{I,S,T\right\}$. Let us define the immediate and expected future reward that the secondary user gains from taking action $i$ as $R^{M}_{i}$ ($M$ for myopic) and $R^{L}_{i}$ ($L$ for long term), respectively, where $i\in\left\{I,S,T\right\}$. Let $p\left(t\right)$ denote the secondary belief that the primary user is idle, $p\left(t\right)\in\left[0,1\right]$. After performing an action $i$ and obtaining a corresponding observation $O$, the updated belief is $\mathcal{E}^{O}_{i}\left(p\right)$. The observations are the sensing outcome or the acknowledgment received from the secondary receiver in case of transmission. If there is no observation associated with the action such as when the secondary user remains idle, the update is $\mathcal{E}_{i}\left(p\right)$. The probability of observing $O$ associated with action $i$ is given by $w^{O}_{i}$. 

The secondary utility function $U_{s}\left(p\left(t\right),t\right)$ is given by
\begin{equation}
U_{s}\left(p\left(t\right),t\right)=
\max\left\{I\left(p\left(t\right),t\right),S\left(p\left(t\right),t\right),T\left(p\left(t\right),t\right)\right\}
\end{equation}
\noindent where $I\left(p\left(t\right),t\right)$, $S\left(p\left(t\right),t\right)$ and $T\left(p\left(t\right),t\right)$ are the secondary user's maximum expected utilities for taking the action of staying idle, carry out spectrum sensing or make data transmission respectively. Next we discuss how to formulate these three utilities depending on the secondary user action.
\subsubsection{First action ``stay idle"}
Although we assume that the secondary user always has data to transmit, staying idle is sometimes the optimal action. It is better for the secondary transmitter, if the primary user is highly likely to be busy, to stay idle and conserve its energy than to consume its energy in sensing because it is more likely that the sensing outcome for the next few actions would be busy. This assumption is based on our primary traffic model which makes it more likely to sense the channel in the same state as the inter-sensing time diminishes \cite{kim}. When the idle action is chosen, the secondary transmitter conserves its energy but on the other hand this causes a reduction in its throughput. We define the following terms:\\
$T_{I}$: Time of staying idle\\
$K_{I}$: Cost of staying idle per unit time\\
For simplification we assume the time unit =1.\\
Therefore the cost of staying idle for $T_{I}$ units is 
\begin{equation}
C_{I}\left(T_{I}\right)=K_{I} T_{I}
\end{equation}
\noindent Parameter $K_{I}$ has units of rate. This is because, as in \cite{senhua_1}, the utility function is mainly based on the secondary throughput. Cost $K_{I}$ is defined as the secondary rate minus the energy saved per unit time expressed in terms of rate. The immediate reward of the secondary user is given by
\begin{equation}
R^{M}_{I}=-C_{I}\left(T_{I}\right)
\end{equation}

Since the action of the secondary user is to stay idle, there are no observations after the time $T_{I}$. The belief is updated as follows
\begin{equation}
\mathcal{E}_{I}\left(p\right)=p\left(t+T_{I}\right)=p\left(t\right) P_{00}\left(T_{I}\right)+\left(1-p\left(t\right)\right) P_{10}\left(T_{I}\right)
\end{equation}
\noindent where $P_{00}\left(t\right)$ is the probability of the channel to be idle at time instant $t+t'$ if it is idle at time $t'$, whereas $P_{10}\left(t\right)$ is the probability of the channel to be idle at time instant $t+t'$ if it is busy at time $t'$. Probabilities $P_{00}\left(t\right)$ and $P_{10}\left(t\right)$ depend on the on and off distributions and are provided in \cite{kim}. See the Appendix for the derivations of $P_{00}\left(t\right)$ and $P_{10}\left(t\right)$ for the uniform distribution as we use it in our numerical simulations.

Therefore the expected future reward that the secondary user gains after staying idle for time $T_{I}$ can be written as:
\begin{equation}
R^{L}_{I}=U_{s}\left(\mathcal{E}_{I}\left(p\right),t+T_{I}\right)
\end{equation}
\noindent We can then write the secondary user maximum expected utility $I\left(p,t\right)$ for taking the action of staying idle as:
\begin{equation}
I\left(p\left(t\right),t\right)=R^{M}_{I}+\beta \, R^{L}_{I}
\label{I}
\end{equation}
%\begin{equation}
%I\left(p\left(t\right),t\right)=-C_{I}\left(T_{I}\right)+\beta \: U_{s}\left(\mathcal{E}_{I}\left(p\right),t+T_{I}\right)
%\end{equation}
\noindent where $\beta\in\left[0,1\right]$ is the discounting factor. At $\beta =0$, the secondary user only care for the immediate reward and does not take the future into account. The value of $\beta$ is usually very close to one.

\subsubsection{Second action ``spectrum sensing"}
The secondary user senses the spectrum to detect spectral vacancies. Sensing has its cost expended to detect the presence of a signal and to acquire a sufficient number of samples to yield reliable results. We define the following parameters that identify the sensing cost:\\
$T_{S}$: Sensing time\\
$K_{S}$: Sensing cost/time unit\\
The sensing cost for $T_{S}$ units as a function of the sensing time is as follows:
\begin{equation}
C_{S}\left(T_{S}\right)=K_{S} T_{S}
\label{sensing_cost}
\end{equation}
The immediate reward is then
\begin{equation}
R^{M}_{S}=-C_{S}\left(T_{S}\right)
\end{equation}
The outcome of the sensing process is either free $\left\{O=F\right\}$ or busy $\left\{O=B\right\}$. Spectrum sensing introduces false alarms and mis-detections, which are decreasing functions of the sensing time $T_{S}$. For a target detection probability $P_d$ , the probability of false alarm is related to the target detection probability as follows \cite{sensing_throughput_tradeoff}:
\begin{equation}
P_{fa}\left(T_S\right)= Q\left(\sqrt{2\psi+1} Q^{-1}\left(P_d\right)+\sqrt{T_S f_s}\psi\right)
\label{P_fa}
\end{equation}
\noindent where $Q\left(.\right)$ is the complementary distribution function of the standard Gaussian and $f_s$ is the sampling frequency. The number of samples used for detecting the primary activity is $f_s \, T_{S}$ and $\psi$ is the received signal-to-noise ratio (SNR).

We introduce now a quantity that is important in the construction of the secondary utility function. This quantity is the probability that the primary user remains idle during the secondary user action given that the primary user is idle \cite{senhua_1}. We denote this conditional probability as $q_{i}\left(t\right)$, where $i$ is a possible secondary action at time $t$. Let $X$ be a random variable describing the duration over which the primary user remains inactive. The probability $q_{i}\left(t\right)$ is given by
\begin{equation}
q_{i}\left(t\right)=\Pr\left\{X >t+T_{i}| X >t\right\}
\end{equation}
\begin{equation}
q_{i}\left(t\right)=\frac{\Pr\left\{X >t+T_{i}\right\}}{\Pr\left\{X >t\right\}}
\end{equation}
\begin{equation}
q_{i}\left(t\right)=\frac{1-F_{X}\left(t+T_{i}\right)}{1-F_{X}\left(t\right)}
\label{eq_q}
\end{equation}
where $F_{X}\left(.\right)$ is the cumulative distribution function of the primary user idle time, and $T_{i}$ is the time duration of the secondary's $i$th action. The characteristics of the $q_{i}\left(t\right)$ function vary according to the distribution of the idle period of the primary user. The function $q_{i}\left(t\right)$ is a decreasing function of time for many distributions such as uniform distribution, Gaussian distribution and Rayleigh distribution because as time increases the probability that the primary user return increases. For an exponential distribution, $q_{i}$ is constant due to the memoryless property.

In case the secondary user decides to carry out spectrum sensing, the probability that the sensing outcome is free is 
\begin{equation}
w^{F}_{S}=p\left(t\right) q_{S}\left(t\right) \left(1-P_{fa}\left(T_{S}\right)\right)+\left(1-p\left(t\right) q_{S}\left(t\right)\right)\left(1-P_{d}\right)
\label{case_2_second}
\end{equation}
\noindent The assumption underlying this formula is that the sensing outcome is free only when the primary remains idle over all the sensing duration. This is a valid assumption given the traffic model so long as the sensing duration is small compared to the mean on/off durations, $T_{\text{\scriptsize off}}$ and $T_{\text{\scriptsize on}}$. The belief update can be written as follows using Bayes' rule.
\begin{equation}
\mathcal{E}^{F}_{S}\left(p\right)=p\left(t+T_{S}\right)=\frac{p\left(t\right) q_{S}\left(t\right) \left(1-P_{fa}\left(T_{S}\right)\right)}{w^{F}_{S}}
\label{case_2_first}
\end{equation}
Similarly, if spectrum sensing is carried out, the probability of a busy outcome is
\begin{equation}
w^{B}_{S}=p\left(t\right) q_{S}\left(t\right) P_{fa}\left(T_{S}\right)+\left(1-p\left(t\right) q_{S}\left(t\right)\right)P_{d}
\label{case_3_second}
\end{equation}
\noindent The belief update is consequently given by
\begin{equation}
\mathcal{E}^{B}_{S}\left(p\right)=p\left(t+T_{S}\right)=\frac{p\left(t\right) q_{S}\left(t\right) P_{fa}\left(T_{S}\right)}{w^{B}_{S}}
\label{case_3_first}
\end{equation}
\noindent The expected future reward when spectrum sensing is carried out
\begin{equation}
R^{L}_{S}=w^{F}_{S} U_{s}\left(\mathcal{E}^F_{S}\left(p\right),t+T_{S}\right)+w^{B}_{S} U_{s}\left(\mathcal{E}_{S}^{B}\left(p\right),t+T_{S}\right)
\end{equation}
And hence we can write the maximum expected utility $S\left(p,t\right)$ that the secondary user gains from sensing the spectrum as:
\begin{equation}
S\left(p(t),t\right)=R^{M}_{S}+\beta \, R^{L}_{S}
\label{S}
\end{equation}
%\begin{equation}
%\begin{array}{rl}
%S\left(p,t\right)~&=~-C_{S}\left(T_{S}\right)+\beta \: w^{F}_{S} \:  %U_{s}\left(\mathcal{E}_{B}^{S}\left(p\right),t+T_{S}\right)\\
%&+~ \beta \: w^{B}_{S} \:  U_{s}\left(\mathcal{E}_{S}^{B}\left(p\right),t+T_{S}\right)
%\end{array}
%\end{equation}
\subsubsection{Third action ``data transmission"}
The immediate reward $R^{M}_{T}$ for the secondary transmitter after transmitting its date differs from the previous actions as there will be a reward for successful transmissions (increment in the secondary user throughput) and a collision cost for colliding with the primary user ( penalizing the interference to the primary user). Define the following terms:\\
$R$: Reward/ time unit for successful transmission\\
$\alpha$: Overhead time\\
$C_{C}$: Collision Cost/ time unit\\
$T_{T}$: Transmission time\\
$K_{T}$: Transmission cost/ time unit\\
The transmission energy can be written as:
\begin{equation}
C_{T}\left(T_{T}\right)=K_{T} T_{T}
\end{equation}
The collision cost $C_{C}$ can be written as a function of a factor $\gamma$ which is controlled by the primary user to prevent the secondary user from transmitting a lot without caring for the primary user, where $\gamma\in\left[0,1\right]$.
When $\gamma=0$, the primary user is afforded maximum protection.
\begin{equation}
C_{C}=C_{C max}\left(1-\gamma\right)
\end{equation}
The secondary receiver sends an acknowledgment to the secondary transmitter upon the processing of the received packet. An ACK is sent for correct decoding, while a NACK means that the receiver has failed in decoding the transmitted message. Thus the observations when the secondary terminal transmits are $\left\{O=A\right\}$ or $\left\{O=N\right\}$. Note that receiving a NACK from the secondary receiver does not mean that a collision with the primary user has occurred, because receiving a NACK may result, for instance, from deep channel fades between the secondary transmitter and secondary receiver. On the other hand it is possible that the secondary receiver can successfully decode the secondary transmitter message even when the primary user is transmitting concurrently. Define the following probabilities:
\\$P_{NC}$: probability that the secondary transmitter receives a NACK although no collision with the primary has occurred.\\
$P_{C}$: probability that the secondary transmitter receives a NACK given that a collision happened with the primary user.\\

\noindent The immediate expected reward that the secondary user gains after data transmission is:

%\begin{equation}
%\begin{array}{rl}
%\mbox{\fontsize{7}{7.5}\selectfont $R^{M}_{T}~$}&\mbox{\fontsize{7}{7.5}\selectfont $=~\overbrace{\left[p(t) \, q_{T}(t) %\left(1-P_{NC}\right)+\left(1-p(t) \, q_{T}(t)\right)\left(1-P_{C}\right)\right] R  \left(T_{T}-\alpha\right) %}^{reward~for~transmission~success}$}\\
%&\mbox{\fontsize{7}{7.5}\selectfont $-~\underbrace{\left(1-p(t) \, q_{T}(t)\right) C_{C}  T_{T}}_{collision~penalty}- %\underbrace{C_{T}\left(T_{T}\right)}_{energy~cost}$}
%\end{array}
%\end{equation} 

\begin{equation}
\begin{array}{rl}
R^{M}_{T}&=\left[p(t) \, q_{T}(t) \left(1-P_{NC}\right)+\left(1-p(t) \, q_{T}(t)\right)\left(1-P_{C}\right)\right]\\
&~R  \left(T_{T}-\alpha\right) 
-\left(1-p(t) \, q_{T}(t)\right) C_{C}  T_{T}- C_{T}\left(T_{T}\right)
\end{array}
\end{equation} 

The probability that the secondary user receives an ACK is
\begin{equation}
w^{A}_{T}=p(t) \, q_{T}(t) \left(1-P_{NC}\right)+\left(1-p(t) \, q_{T}(t)\right)\left(1-P_{C}\right)
\end{equation}
\noindent The probability that the primary user is idle after $T_{T}$ given that the secondary transmitter receives an ACK is:
\begin{equation}
\mathcal{E}_{T}^{A}\left(p\right)=p\left(t+T_{T}\right)=\frac{p(t) \, q_{T}(t) \left(1-P_{NC}\right)}{w^{A}_{T}}
\end{equation}
Similarly, the probability that the secondary user receives a NACK is
\begin{equation}
w^{N}_{T}=p(t) \, q_{T}(t) \, P_{NC}+\left(1-p(t) \, q_{T}(t)\right)P_{C}
\end{equation}
\noindent The update in the NACK case is
\begin{equation}
\mathcal{E}^{N}_{T}\left(p\right)=p\left(t+T_{T}\right)=\frac{p(t) \,q_{T}(t) P_{NC}}{w^{N}_{T}}
\end{equation}
The expected future reward is
\begin{equation}
R^{L}_{T}=w^{A}_{T} \, U_{s}\left(\mathcal{E}_{T}^{A}\left(p\right),t+T_{T}\right)+w^{N}_{T} \, U_{s}\left(\mathcal{E}^{N}_{T}\left(p\right),t+T_{T}\right)
\end{equation}
\noindent Finally,
\begin{equation}
T\left(p(t),t\right)=R^{M}_{T}+\beta \, R^{L}_{T}
\label{T}
\end{equation}
%\begin{equation}
%\begin{array}{rl}
%T\left(p,t\right)~&=~\left(p \, q_{T} \left(1-P_{NC}\right)+\left(1-p \, q_{T}\right)\left(1-P_{C}\right)\right)\\
%&~ R  \left(T_{T}-\alpha\right)-\left(1-p \, q_{T}\right) C_{C}  T_{T}\\ 
%&-~ C_{T}\left(T_{T}\right) + \beta \, w_{T}^{A} \: U_{s}\left(\mathcal{E}^{A}_{T}\left(p\right),t+T_{T}\right)\\
%&+~ \beta \, w_{T}^{N} \: U_{s}\left(\mathcal{E}^{N}_{T}\left(p\right),t+T_{T}\right)
%\end{array}
%\end{equation}
%------------------------------------------------Problem formulation------------------------------------------------------------
\section{Adaptive Sensing and Transmission}{\label{section:adaptive_structure}}
\begin{figure}[tbp]
	\centering
		\includegraphics[width=0.4\textwidth]{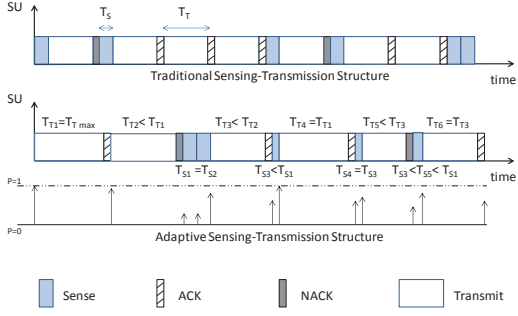}
	\caption{Traditional vs. Adaptive Sensing-Transmission Structures}
	\label{fig:idea}
\end{figure}
Many existing works on cognitive radio networks assume that the secondary users have fixed sensing and transmission durations as illustrated in Fig. \ref{fig:idea}. This is not necessarily optimal as the secondary transmitter should not waste time and energy in long sensing periods although it has a high belief that the primary user is idle. This time can be exploited in relatively long transmissions, thereby increasing secondary throughput and maximizing its utility. On the other hand, at low values of $p$ the secondary transmitter can increase its sensing time to get reliable results on the primary user state as increasing the sensing time decreases the false alarm probability and the probability of mis-detection. In this case, it can also decrease its transmission duration to reduce the probability of collision with the primary user.

The secondary transmitter can adaptively vary its sensing and transmission durations according to its belief about the primary user state $p\left(t\right)$. In the sequel, we fix $T_{I}$ and assume that the sensing and transmission times are linear functions of $p$. That is,
\begin{equation}
T_{T}\left(p\right)=a_{0}+a_{1}p
\end{equation}
\begin{equation}
T_{S}\left(p\right)=b_{0}-b_{1}p
\end{equation}
\noindent where $a_{0}$, $a_{1}$, $b_{0}$ and $b_{1}$ are our design parameters. Note that if a fixed $T_{T}$ or $T_{S}$ is optimal, then the solution of the optimization problem would yield close-to-zero $a_{1}$ or $b_{1}$.

Parameters $a_{0}$, $a_{1}$, $b_{0}$ and $b_{1}$ are nonnegative and obey the following inequalities:
%\begin{eqnarray}
%\label{parameter_constraints} 
% \quad && T_{T,min} \leq a_{0} \leq a_{0,max} \nonumber  \\
%      && T_{S,min} \leq b_{0} \leq T_{S,max}  \\
%      && 0 \leq a_{1} \leq T_{T,max} - a_{0,max} \nonumber \\
%      && -b_{0} \leq b_{1} \leq 0 \nonumber
%\end{eqnarray}
\begin{eqnarray}
\label{parameter_constraints} 
 \quad &&  a_{0} \geq T_{T,{\rm min}} \nonumber  \\
      &&  b_{0}-b_{1} \geq T_{S,{\rm min}}  \\
      && a_{0}+ a_{1} \leq T_{T,{\rm max}}  \nonumber \\
      && b_{0} \leq T_{S,{\rm max}} \nonumber
\end{eqnarray}
\noindent where the positive parameters $T_{T,{\rm min}}$ , $T_{S,{\rm min}}$ , $T_{T,{\rm max}}$ and $T_{S,{\rm max}}$ are the minimum and maximum transmission and sensing durations, respectively. The minimum duration for transmission is dictated by data overhead and the shortest possible data payload, whereas the minimum time for sensing is dictated by a minimal detection reliability requirement. The maximum durations are imposed to protect the primary by frequently checking its activity. Moreover, the maximum durations should be considerably less than $T_{\text{\scriptsize off}}$ and $T_{\text{\scriptsize on}}$ in order for the probability formulas to be valid. The first two constraints in (\ref{parameter_constraints}) maintain that the sensing and transmission durations are not less than the minimum specified values for all $p\left(t\right) \in [0,1]$. 

Our objective is to obtain the optimal sensing and transmission durations to maximize the secondary user utility. In other words, our objective is to dynamically decide, for each secondary user action, the optimal sensing and transmission durations to maximize the secondary user utility. Define
\begin{equation}
\begin{split}
&S^{*}\left(p\left(t\right),t\right)=\max_{a_{0},a_{1},b_{0},b_{1}} S\left(p\left(t\right),t\right)\\
&T^{*}\left(p\left(t\right),t\right)=\max_{a_{0},a_{1},b_{0},b_{1}} T\left(p\left(t\right),t\right)
\end{split}
\end{equation}
\noindent In the adaptive case, the secondary utility then becomes
\begin{equation}
\begin{split}
&U_{s}\left(p\left(t\right),t\right)=\\
&\max \left\{ I\left(p\left(t\right),t\right), S^{*}\left(p\left(t\right),t\right), T^{*}\left(p\left(t\right),t\right)\right\}
\label{U_function}
\end{split}
\end{equation}
\noindent The optimal parameters if the secondary action is to sense or to transmit are $a^{*}_{0},a^{*}_{1},b^{*}_{0},b^{*}_{1}$ such that
\begin{equation}
\begin{split}
&a^{*}_{0},a^{*}_{1},b^{*}_{0},b^{*}_{1}=\underset{a_{0},a_{1},b_{0},b_{1}}{\mbox{argmax}} S\left(p\left(t\right),t\right)\nonumber \\
& \mbox{if } S^{*}\left(p\left(t\right),t\right) > \max \left\{I\left(p\left(t\right),t\right), T^{*}\left(p\left(t\right),t\right)\right\} \nonumber \\
&\mbox{or }a^{*}_{0},a^{*}_{1},b^{*}_{0},b^{*}_{1}=\underset{a_{0},a_{1},b_{0},b_{1}}{\mbox{argmax}} T\left(p\left(t\right),t\right)\nonumber \\
& \mbox{if } T^{*}\left(p\left(t\right),t\right) > \max \left\{I\left(p\left(t\right),t\right), S^{*}\left(p\left(t\right),t\right)\right\} \nonumber
\end{split}
\end{equation}

The optimal action for the secondary user can be found either by using value iteration or backward induction. Employing value iteration with $\beta < 1$, we initialize the matrix $U_s\left(p\left(t\right),t\right)$ with zeros. We iterate using (\ref{U_function}) until convergence \cite{Bertsekas_DP}. Hence we obtain the optimal action for each $p$ and $t$. Another method which is used in this paper for the numerical results is to use backward induction with $\beta=1$. The possibility of doing backward induction with a unity discounting factor is predicated on the monotonically decreasing nature of the function $q_i\left(t\right)$ given by (\ref{eq_q}) for some distributions. That is, as time proceeds, the probability of the primary user remaining idle during the sensing or transmission phase approaches zero. This means that regardless of the value of $p$, at large times, the secondary utility function is given by
\begin{equation}
U_s\left(p,t\right)=\max\left\{R^{M}_{I},R^{M}_{S},R^{M}_{T}\right\}
\end{equation}
\noindent for large $t$ such that $q_i\left(t\right)$ is almost zero. Given these values, backward induction can be used to get all $U_s$ values at different $p$ and $t$.

As shown in the next section, our solution is a threshold-based policy as in \cite{senhua_1}. This means that the secondary transmits when its belief about the primary being idle exceeds a certain threshold.

%--------------------------------------Fixed sensing and transmission structure-------------------------------------------------
\begin{figure}[tbp]
	\centering
		\includegraphics[width=0.4\textwidth]{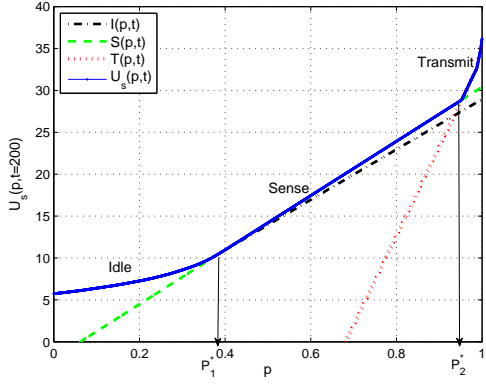}
	\caption{Secondary user utility $U_s(p,t)$ as a function of $p$ at $t=200$. The three components of the utility function, $I(p,t)$, $S(p,t)$, and $T(p,t)$ are depicted. The figure shows the threshold-based nature of the optimal policy.}
	\label{fig:utility}
\end{figure}
\begin{figure}[tbp]
	\centering
		\includegraphics[width=0.4\textwidth]{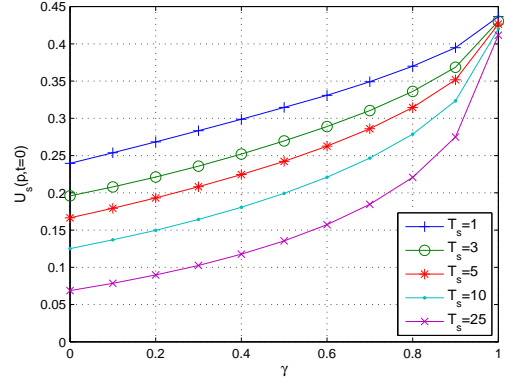}
	\caption{Secondary utility $U_s(p,t)$ as a function of $\gamma$ with $T_T=10$, $T_I=5$, and various values for $T_S$. The case here is non-adaptive transmission and sensing durations with perfect sensing and no data overhead. Utility $U_s(p,t)$ increases with $\gamma$ as the collision penalty decreases. A higher sensing duration is a waste of transmission opportunities as sensing is assumed to be perfect.}
	\label{fig:trad_perfect_1_no_overhead}
\end{figure}
\begin{figure}[tbp]
	\centering
		\includegraphics[width=0.4\textwidth]{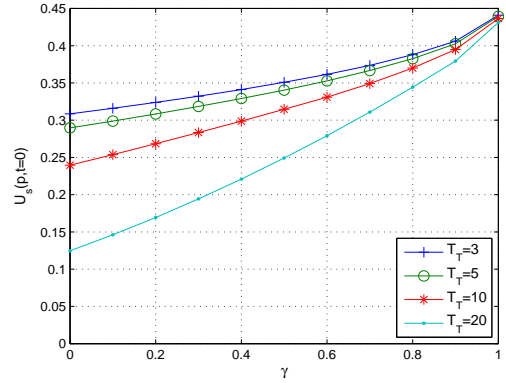}
	\caption{Secondary utility $U_s(p,t)$ as a function of $\gamma$ with $T_S=1$, $T_I=5$, and various values for $T_T$. The case here is non-adaptive transmission and sensing durations with perfect sensing and no data overhead. Utility $U_s(p,t)$ increases with $\gamma$ as the collision penalty decreases. A higher transmission duration reduces the secondary utility due to the increase in collision probability. Hence, the reduction in $U_s(p,t)$ caused by a high $T_T$ decreases as $\gamma$ increases.}
	\label{fig:trad_perfect_2_no_overhead}
\end{figure}

\begin{comment}
$I\left(p,t\right)$ and $T\left(p,t\right)$ are described by equations $\left(\ref{I}\right)$ and $\left(\ref{T}\right)$ respectively. \\
\end{comment}
%------------------------------------Adaptive sensing and transmission structure------------------------------------------------
%\begin{figure}[tbp]
%	\centering
%		\includegraphics[width=0.4\textwidth]{trad_perfect_1.eps}
%	\caption{Effect of varying $T_S$ at fixed $T_T=10$ and $T_I=5$ in the traditional structure, perfect sensing}
%	\label{fig:trad_perfect_1}
%\end{figure}

\begin{figure}[tbp]
	\centering
		\includegraphics[width=0.4\textwidth]{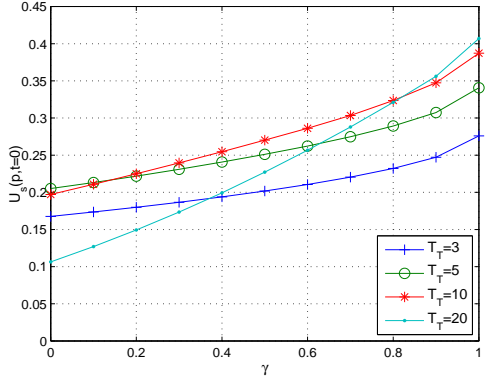}
	\caption{Secondary utility versus $\gamma$ with $T_S=1$, $T_I=5$, and various $T_T$ values under the non-adaptive scheme assuming perfect sensing, but considering data overhead. At small $\gamma$ values, long transmission durations reduce utility due to collision penalty. As $\gamma$ increases, the collision penalty decreases and short transmission durations result in a low utility due to the overhead.}
	\label{fig:trad_perfect_2}
\end{figure}
%----------------------------------------------Simulation results---------------------------------------------------------------
\begin{figure}[tbp]
	\centering
		\includegraphics[width=0.4\textwidth]{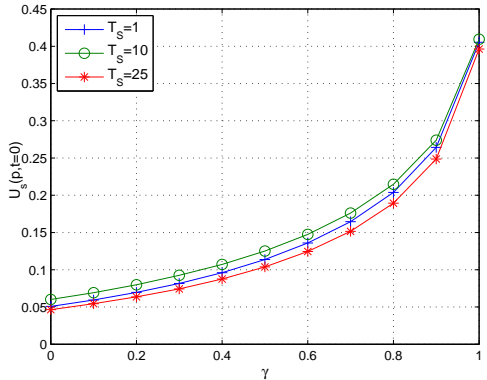}
	\caption{This figure shows the same situation as Fig. \ref{fig:trad_perfect_1_no_overhead} but allowing for sensing errors and assuming a data overhead. Note that the best performance corresponds to the intermediate value for $T_S$. The reason is that although increasing the sensing duration comes at the expense of transmission duration, it increases decision reliability and reduces lost transmission opportunities caused by false alarm.}
	\label{fig:trad_imperfect_1}
\end{figure}
\section{Simulation Results}{\label{section:simulation_results}}
The simulation results section has three main parts. First we show the relation between the secondary user utility $U_s\left(p,t\right)$ and the belief state $p$ at certain time $t$. In the second part, we provide simulation results for the traditional scheme where the sensing and transmission durations for the secondary transmitter are fixed. Finally, we simulate our adaptive scheme and compare it with the fixed one. In the second and third part we study perfect/imperfect sensing with/without the overhead. The results demonstrate the performance enhancement due to adapting the sensing and transmission durations. 

For the results below in the three parts we use the following simulation parameters unless otherwise mentioned. For the idle action, the secondary user idle duration is $T_I=5$ , cost of staying idle per time unit $K_I=0.001$. For the sensing action, we use the sensing cost per time unit $K_S=0.1$. The transmission action parameters are as follows: transmission cost per time unit $K_T=0.1$, overhead time $\alpha=1$, reward per time unit for successful transmission $R=1$. As mentioned before, $0 \leq\gamma\leq 1$ is a factor controlled by the primary user to control the secondary user transmissions by changing the collision cost. At $\gamma=0$, a maximum protection is required which is equivalent to $C_C=C_{C max}=20$, the other extreme case happened at $\gamma=1$ which gives $C_C=0$. we assume that $P_{NC}=0$ and $P_C=1$ for all simulation parts. For the primary user, $f_{\rm on}\left(t\right)$ and $f_{\rm off}\left(t\right)$ are uniform over the interval $\left[0,1000\right]$. We use also backward induction where $\beta=1$. 

\subsection{Secondary user utility function characteristics}
Fig. \ref{fig:utility} shows the relation between the secondary user utility $U_s\left(p,t=200\right)$ and the belief state $p$ in the perfect sensing with overhead case. Note that $U_s\left(p,t=200\right)$ is a convex and increasing function in $p$. It can be shown following an argument similar to that in \cite{senhua_1} that $U_s\left(p,t\right)$ is a convex and increasing function in $p$ at any value of $t$. For this figure, we use a fixed sensing time $T_S=20$ and a fixed transmission time $T_T=7$. For the collision cost we set $\gamma=0.5$.

The threshold based structure is obvious in Fig. \ref{fig:utility} where $p^*_1\left(t=200\right)=0.3939$ and $p^*_2\left(t=200\right)=0.9522$. At $p<p^*_1$, $I\left(p,t\right)$ is greater than $S\left(p,t\right)$ and $T\left(p,t\right)$ which means that the optimal action for the secondary user is to stay idle. At $p^*_1\leq p\leq p^*_2$, $S\left(p,t\right)$ is greater than $I\left(p,t\right)$ and $T\left(p,t\right)$ which means that the optimal action for the secondary user is to sense the spectrum. At $p>p^*_2$, the optimal action is to transmit. 

\subsection{Traditional scheme} Here we consider the case where the secondary transmitter has a fixed sensing and transmission durations for all its access period. We show the variation of the secondary user utility with $\gamma$ at different values for the fixed sensing and transmission durations. In this part we consider perfect/imperfect sensing with/without the overhead.  

\subsubsection{Perfect sensing without overhead} The change of the secondary user utility with $\gamma$ at different sensing durations is shown in Fig. \ref{fig:trad_perfect_1_no_overhead}. We fix the transmission duration at $T_T=10$. Increasing the sensing duration for the secondary transmitter decreases its utility since we assume in this part that the secondary user has a perfect sensing mechanism. That is, the optimal sensing time for the secondary user is $T_S=T_{S min}=1$. The secondary user obtains no gain from increasing the sensing duration as $P_{fa}=0$ and $P_d=1$ and, in fact, sensing wastes time that can potentially be used for data transmission. We notice from the figure that the secondary user utility $U_s\left(p,t\right)$ increases with $\gamma$, because as $\gamma$ increases, the collision cost decreases and, hence, the utility increases.

The effect of varying the transmission duration on the $U_s$ versus $\gamma$ curve at a fixed sensing time $T_S=1$ is shown in Fig. \ref{fig:trad_perfect_2_no_overhead}. When $\alpha=0$, the secondary user utility decreases as the transmission duration increases at fixed sensing duration. This decrease is reduced as $\gamma$ increases. The reason for this is that a longer transmission duration means a higher probability of colliding with primary transmission. Since when $\gamma$ increases, the collision penalty decreases, the degradation caused by a long transmission duration is reduced.

\subsubsection{Perfect sensing with overhead} We now study the secondary utility considering the data overhead. The case corresponding to Fig. \ref{fig:trad_perfect_1_no_overhead} with overhead is not provided here as the utility is reduced at all $\gamma$ due to overhead. However the effect of varying the transmission duration at a fixed sensing time $T_S=1$ produces an interesting change as shown in Fig. \ref{fig:trad_perfect_2}. When the collision penalty is low at high $\gamma$'s, small transmission durations barely above the data overhead result in small secondary utility relative to high transmission durations. On the other hand, at high collision costs, long transmission durations reduce the secondary utility. 

\begin{figure}[tbp]
	\centering
		\includegraphics[width=0.4\textwidth]{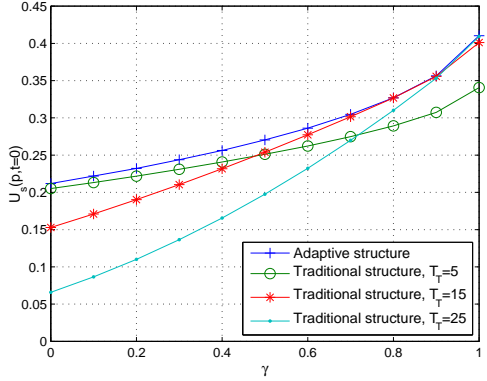}
	\caption{Comparing the adaptive scheme with the traditional one at fixed $T_I=5$ and varying $T_T$ and $T_S$, assuming perfect sensing with data overhead. Note that the improvement in the secondary utility is due to the adaptive transmission duration not the varying sensing duration.}
	\label{fig:adap_perfect_1}
\end{figure}

\begin{figure}[tbp]
	\centering
		\includegraphics[width=0.4\textwidth]{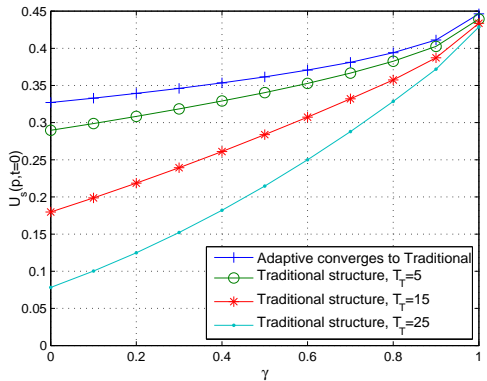}
	\caption{Comparing the adaptive scheme with the traditional one at fixed $T_I=5$ and varying $T_T$ and $T_S$ under perfect sensing without data overhead. In this case only, the adaptive structure converges to the non-adaptive one where the optimal sensing and transmission durations are the fixed ones.}
	\label{fig:adap_perfect_1_no_overhead}
\end{figure}
%-------------------------------------------------------------------------------------------------------------------------------
\begin{figure}[tbp]
	\centering
		\includegraphics[width=0.4\textwidth]{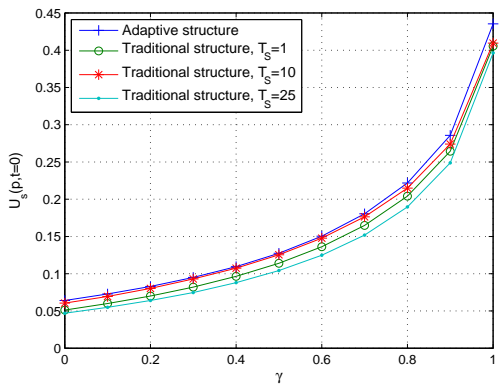}
	\caption{Comparing the adaptive scheme with the traditional one at fixed $T_I=5$ and varying $T_T$ and $T_S$. Sensing is imperfect with overhead. Note that in this case, the optimal durations are the adaptive ones.}
	\label{fig:adap_imperfect_1}
\end{figure}
\subsubsection{Imperfect sensing with overhead} Recall that the false alarm probability varies with the sensing time according to (\ref{P_fa}). We use the following parameters in the imperfect sensing part: detection probability $P_d=0.9$, $B.W.=50 KHz$, sampling frequency $f_s=5/8 B.W.$, $\psi=-25 dB$ and $K_S=0.01$. Fig. \ref{fig:trad_imperfect_1} shows the performance of the secondary user utility while varying the sensing duration at fixed transmission time $T_T=10$. The case here is different from Fig. \ref{fig:trad_perfect_1_no_overhead}. As previously mentioned, under perfect sensing, the optimal sensing duration is the minimum possible. Increasing the sensing duration brings no additional benefit or refinement of the sensing outcome. In the imperfect sensing case, there is a tradeoff between increasing the sensing time and consequently getting a lower probability of false alarm, and decreasing it to have more time for transmission.

This tradeoff is well demonstrated in Fig. \ref{fig:trad_imperfect_1}. As the sensing duration is increased, we get a higher utility due to the lower false alarm probability. However an excessive increase in the sensing duration degrades the performance because little time is left for transmission. 

The effect of varying the transmission duration for the secondary user while fixing the sensing duration is not depicted here as it is similar to Fig. \ref{fig:trad_perfect_2} with a degradation of utility at all collision costs due to imperfect sensing.
Simulations for the imperfect sensing without the overhead are omitted.

\subsection{Adaptive scheme} In this part we compare our adaptive sensing and transmission scheme with the traditional fixed one. We show that using adaptive durations for sensing and transmission returns a higher utility for the secondary user. We set $T_{T,{\rm min}}=1$, $T_{T,{\rm max}}=30$, $T_{s,{\rm min}}=1$, $T_{s,{\rm max}}=10$.

\subsubsection{Perfect sensing with overhead} The effect of varying the sensing and transmission durations for the secondary user at each value of $\gamma$ is shown in Fig. \ref{fig:adap_perfect_1} where we fix the sensing duration in the traditional structure at $T_S=1$ and simulate the system at different transmission durations. We notice that at every value of $\gamma$ which corresponds to a certain collision cost, the adaptive structure returns a higher utility for the secondary user than using fixed durations. However solving our constraint optimization problem in the perfect sensing case gives us a fixed sensing duration and an adaptive transmission one. That is the optimizer returns $b_0^*=1$ and $b_1^*=0$ for all values of $\gamma$ which gives a sensing duration as:
\begin{equation}
T_S=b_0^*+b_1^* p=1
\label{min_sensing}
\end{equation}
This is what we expect for perfect sensing case.

Regarding transmission durations, solving the optimization problem returns different values for $a_0^*$ and $a_1^*$ for every value of $\gamma$. The higher utility in Fig. \ref{fig:adap_perfect_1} is due to the adaptive transmission duration at every value of $\gamma$, not the varying sensing durations. This is not the case under imperfect sensing where the optimizer chooses $b_1^*\neq0$.

\subsubsection{Perfect sensing without overhead} Here $b_0^*=1$ and $b_1^*=0$ for all values of $\gamma$ which gives a sensing duration as in equation (\ref{min_sensing}). Also for the transmission time, the optimizer returns $a_0^*=1$ and $a_1^*=0$ for all values of $\gamma$ which gives the minimum transmission duration as expected due to the absence of the overhead.  
\begin{equation}
T_T=a_0^*+a_1^* p=1
\end{equation}
Actually in this case only, the adaptive structure converges to the fixed one as shown in Fig. \ref{fig:adap_perfect_1_no_overhead}. i.e. the secondary transmitter uses fixed sensing and transmission durations in order to maximize its utility.

\subsubsection{Imperfect sensing with overhead} To test the impact of adapting both the sensing and transmission durations for the secondary user we consider in this part the imperfect sensing with overhead case. Solving our optimization problem, we found that both the sensing and the transmission durations are varying adaptively according to the belief that the primary user is idle and that $a_0^*$, $a_1^*$, $b_0^*$ and $b_1^*$ will all have values that do not equal to zero at all values of $\gamma$. Fig. \ref{fig:adap_imperfect_1} shows that it is better for the secondary transmitter to adaptively change its sensing and transmission durations according to the belief $p$ as that increases secondary utility. 
We notice that at $\gamma=1$, the optimizer chooses the maximum transmission time and the minimum sensing time as follows because at this value of $\gamma$ there is no collision penalty.
\begin{equation}
T_T=T_{T max}=a_0^*+a_1^*p
\end{equation}
\begin{equation}
T_S=T_{S min}=b_0^*-b_1^*p
\end{equation}
%---------------------------------------------------Conclusion------------------------------------------------------------------
\section{Conclusion}{\label{section:conclusion}}
We have developed an adaptive scheme for the sensing and transmission durations of a secondary user sharing a channel with a primary user. The sensing and transmission durations are varied adaptively according to the secondary belief regarding primary activity. The objective is to maximize the secondary utility which takes into account the impact of the secondary user's decision on the future. Simulation results have demonstrated that the proposed adaptive scheme returns a higher utility than the non-adaptive one.

Several interesting directions for future work exist. For example, a power control scheme can be incorporated so that the secondary adapts its transmission power based on its belief regarding the primary state of activity and also the channels connecting the primary and secondary transmitters and receivers. Furthermore we can consider the case where there is a kind of cooperation between primary and secondary users. The incentive for the primary user would be some extra revenue or some help from the secondary user in relaying its message. Finally, the investigation can be made more realistic by incorporating the primary and secondary queues in the analysis.

%---------------------------------------------------Appendix--------------------------------------------------------------------
\appendix 
\subsection{Derivation of $P_{00}\left(t\right)$ and $P_{10}\left(t\right)$}
Probabilities $P_{00}\left(t\right)$ and $P_{10}\left(t\right)$ depend on the on and off distributions of the primary user. Using renewal theory, $P_{11}\left(t\right)$ can be expressed as:
\begin{equation}
P_{11}\left(t\right)=\int_{t}^{\infty} \frac{f_{\rm on}\left(u\right)}{T_{\rm on}} \,du + \int_{0}^{t} h_{10}\left(u\right) f_{\rm on}\left(t-u\right) \,du
\label{P_11_t}
\end{equation}
where $h_{10}\left(u\right)$ is the renewal density of the off state given that the renewal process started from the on state. It is proven in \cite{cox_renewal} that $h_{10}\left(s\right)$ is given by: 
\begin{equation}
h_{10}\left(s\right)=\frac{f_{\rm off}\left(s\right) \left(1-f_{\rm on}\left(s\right)\right)} {T_{\rm on} s \left(1-f_{\rm on}\left(s\right) f_{\rm off}\left(s\right) \right)}
\end{equation}
By applying Laplace transform to equation $\left(\ref{P_11_t}\right)$, we get
\begin{equation}
P_{11}\left(s\right)=\frac{1}{s} -\frac{\left(1-f_{\rm on}\left(s\right)\right) \left(1-f_{\rm off}\left(s\right)\right)} {T_{\rm on} s^2 \left(1-f_{\rm on}\left(s\right) f_{\rm off}\left(s\right) \right)}
\label{P_11_s}
\end{equation}
Our objective now is to derive the formulas of $P_{10}\left(t\right)$ and $P_{00}\left(t\right)$. Using the inverse Laplace transform of equation $\left(\ref{P_11_s}\right)$, we can get $P_{10}\left(t\right)$ as:
\begin{equation}
P_{10}\left(t\right)=1-P_{11}\left(t\right)
\end{equation}
Similarly, using the inverse Laplace transform of equation $\left(\ref{P_00_s}\right)$ we can get $P_{00}\left(t\right)$ as
\begin{equation}
P_{00}\left(s\right)=\frac{1}{s} -\frac{\left(1-f_{\rm off}\left(s\right)\right) \left(1-f_{\rm on}\left(s\right)\right)} {T_{\rm off} s^2 \left(1-f_{\rm off}\left(s\right) f_{\rm on}\left(s\right) \right)}
\label{P_00_s}
\end{equation}
We now focus on the case when $f_{\rm on}\left(t\right)$ and $f_{\rm off}\left(t\right)$ are uniformly distributed on the interval $\left[0,b\right]$
\begin{equation}
T_{\rm on}=T_{\rm off}=\frac{b}{2}\nonumber\\
\end{equation}
\begin{equation}
f_{\rm on}\left(t\right)=f_{\rm off}\left(t\right)=\frac{1}{b} \left[u\left(t\right)-u\left(t-b\right)\right]
\end{equation}
Using Laplace transform we can get:
\begin{equation}
f\left(s\right)=f_{\rm on}\left(s\right)=f_{\rm off}\left(s\right)=\frac{1}{bs} \left[1-\exp\left(-bs\right)\right]
\end{equation}
\begin{equation}
P_{11}\left(s\right)=\frac{1}{s} -\frac{\left(1-f\left(s\right)\right)^2 } {T_{\rm on} s^2 \left(1-f^2\left(s\right) \right)}\nonumber\\
\end{equation}
\begin{equation}
=\frac{1}{s} -\frac{\left(1-f\left(s\right)\right) } {T_{\rm on} s^2 \left(1+f\left(s\right) \right)}\nonumber\\
\end{equation}
\begin{equation}
=\frac{1}{s} -\frac{2} {T_{\rm on} s^2 \left(1+f\left(s\right) \right)}- \frac{1}{T_{\rm on} s^2}\nonumber\\
\end{equation}
\begin{equation}
=\frac{1}{s} -\frac{4} {s\left(bs+1-\exp\left(-bs\right)\right)}- \frac{1}{\left(b/2\right) s^2}
\end{equation}
Using the Maclaurin series for $\frac{1}{1-x}=1+x+x^2+x^3+...$ , where $x=\frac{\exp\left(-bs\right)}{1+bs}$ we get the following:
\begin{equation}
\begin{array}{rl}
P_{11}\left(s\right)~&=~\frac{1}{s} -\frac{1}{\left(b/2\right) s^2} -(\frac{4} {s\left(1+bs\right)}+\frac{4 \exp\left(-bs\right)} {s\left(1+bs\right)^2} \\
&+~\frac{4 \exp\left(-2bs\right)} {s\left(1+bs\right)^3} +\frac{4 \exp\left(-3bs\right)} {s\left(1+bs\right)^4}+...) \nonumber\\
\end{array}
\end{equation}
\begin{equation}
\begin{array}{rl}
~&=~\frac{1}{s} -\frac{1}{\left(b/2\right) s^2} - (g_0(s) + g_1(s)\exp(-bs) \\
&+~ g_2(s)\exp(-2bs) + g_3(s)\exp(-3bs) + ...)
\end{array}
\end{equation}
After some algebraic computations we get the following: 
\begin{equation}
\begin{array}{rl}
P_{11}(t)~&=~u(t)-\frac{2t}{b}u(t)-(g_0(t)+g_1(t-b)u(t-b)\\
&+~g_2(t-2b)u(t-2b)+g_3(t-3b)u(t-3b)+...)
\end{array}
\end{equation}
since
\begin{equation}
1<t<b\nonumber\\
\end{equation}
Therefore, we get only the first term of the series $g_0(t)$
\begin{equation}
P_{11}(t)=u(t)-\frac{2t}{b}u(t)-(4-4 \exp(-t/b))
\end{equation}

\begin{equation}
P_{10}\left(t\right)=1-P_{11}\left(t\right)\nonumber\\
\end{equation}
then at $t=T_I$ we get the following:
\begin{equation}
P_{10}(T_I)= (5-4 \exp(-T_I/b))-u(T_I)+\frac{2T_I}{b}u(T_I) 
\end{equation}
Following the same argument for $P_{00}(s)$, we can get $P_{00}(T_I)$ as follows:
\begin{equation}
P_{00}(T_I)=u(T_I)-\frac{2T_I}{b}u(T_I)-(4-4 \exp(-T_I/b))
\end{equation}
%---------------------------------------------------References------------------------------------------------------------------
\bibliographystyle{IEEEtran}
\bibliography{IEEEabrv,MyLib}

\end{document}